\documentclass[11pt,hyp,]{amsart}
\usepackage{amsmath,amsfonts,a4,epsfig}
\usepackage{graphicx}
\usepackage[latin1]{inputenc}
\usepackage[french, english]{babel}
\usepackage{color}
\usepackage{enumerate}
\usepackage{setspace}
\usepackage{fancyhdr}

\def\preuve{\begin{proof}}

\newtheorem{defi}{Definition}[section]
\newtheorem{defini}{Definitions}[section]

\newtheorem{lemm}{Lemma}[section]
\newtheorem{prop}{Proposition}[section]
\newtheorem{rem}{Remark}[section]

\newtheorem{coro}{Corollary}[section]
\newtheorem{theo}{Theorem}[section]

  \newenvironment{demo}{\noindent {\it Proof:}
      \begin{quotation}\noindent}{\end{quotation}\hfill$\square $}

\usepackage{hyperref}
\hypersetup{nesting=true,debug=true,naturalnames=true}
\usepackage{graphicx,amssymb,upref}

\let\<\langle
\let\>\rangle
\usepackage[all,pdf]{xy}
\UseComputerModernTips

\let\uml\"

\title[Sectoriality and essential spectrum]{Sectoriality and essential spectrum of non symmetric graph Laplacians }

\author{MARWA BALTI, COLETTE ANN\'E, NABILA TORKI-HAMZA}

\address{Universit\'e de Nantes, Laboratoire de Math\'ematique Jean Lauray, CNRS, Facult\'e des Sciences, BP 92208, 44322 Nantes, (France).}
\email{ colette.anne@univ-nantes.fr}

\address{Universit\'e de Carthage, Facult\'e des Sciences de Bizerte: Math\'ematiques et Applications (UR/13ES47) 7021-Bizerte (Tunisie)\\
Universit\'e de Nantes, Laboratoire de Math\'ematique Jean Lauray, CNRS, Facult\'e des Sciences, BP 92208, 44322 Nantes, (France).
}
\email{ balti-marwa@hotmail.fr}
 
\address{Universit\'e de Kairouan, Isig-K, 3100-Kairouan; Tunisie.}
\email{ natorki@gmail.com}

\keywords{Directed graph, Non self-adjoint Laplacian, Numerical range, Sectorial operator, Essential spectrum.}

\subjclass[2010]{47A45, 47A12, 47A10, 47B25, 39A12, 47B37.}

\begin{document}
\begin{abstract}
We consider a non self-adjoint Laplacian on a directed graph with non symmetric edge weights. We give necessary conditions for this Laplacian to be sectorial. We introduce a special self-adjoint operator and compare its essential spectrum with that of the non self-adjoint Laplacian considered.
\end{abstract}
\maketitle
\tableofcontents
\section*{Introduction}
 Many problems in mathematical physics can be modeled by using operators in Hilbert spaces. Under certain circumstances these operators will be self-adjoint, but in many examples as, e.g., for Schr\"odinger and Sturm-Liouville problems with complex, non-symmetric coefficients and for the discret Laplacian with non symmetric weights, this is not the case. Looking for other approaches, we find that the class of so-called sectorial operators provides an appropriate framework which suits these problems very well. 
 Linear operators of sectorial type play a prominent role in the theory of semigroups and evolution equations; see for instance \cite{Ar}, \cite{Yos}, \cite{Kts}. The sectorial operators are given special attention in view of later applications to the spectral theory and to the analytic and the asymptotic perturbation theory \cite{Kts}, \cite{Ev}, \cite{Kh}. 
Our prime concern is the study of the essential spectrum of a sectorial discrete Laplacian.
We give a criterion for the sectoriality of the graph Laplacian and we compare its essential spectrum with the essential spectrum of its real part. We establish the emptiness of the essential spectrum for both operators. Fujiwara \cite{Fuji}, Keller \cite{Klr} and Balti \cite{Ma} introduced a criterion for the lack of essential spectrum of the graph Laplacian using the positivity of a Cheeger constant at infinity. We extend these results for weighted non symmetric graphs with any non negative Cheeger constant at infinity.

Section \ref{In} is devoted to some definitions and notions on a directed graph with non symmetric edge weights
and the associated non symmetric differential Laplacian $\Delta$. We recall some basic results introduced in \cite{Ma} like Green's formula and the spectral properties of $\Delta$ and of its formal adjoint for weighted graphs satisfying the Asumption $(\beta)$ which assure to conclude an equality of in and out conductivity at each vertex.\\
 In section \ref{2} we state that under a geometric hypothesis $(\gamma)$ see page 5, the closure operator $\overline{\Delta}$  is sectorial. We use the benefit of this notion to characterize the essential spectrum.\\
  In section \ref{sym}, we propose a general condition, to ensure self-adjointness
for the symmetrized operator $\overline{\Delta}+\Delta^{*}$. Then, we use the comparison Theorem of Lewis to show the relation between the essential spectra of $\overline{\Delta}$ and $\overline{\Delta}+\Delta^{*}$ and the absence of the essential spectrum for both.\\
In section \ref{EX}, we show the importance of the assumption $(\gamma)$ for the emptiness of the essential spectrum by the study of a counter-example. Then, we propose an other example of graph such that the essential spectrum is empty in spite of the Cheeger constant at infinity is zero.
\section{Preliminaries}\label{In}
We review in this section some basic definitions which are developed in \cite{Ma}, \cite{Mar} and we introduce the notations used in the article.
\subsection{Notion of Graphs}
A directed weighted graph is a couple$\left( G,b\right) $ and $G=(V,\vec{E})$, where $V$ is a countable set of vertices, $\vec{E}\subset V\times V$ is the set of directed edges and $b : V \times V\to\left[0,\infty\right)  $ is a weight satisfying the following conditions:
\begin{itemize}
  \item $b(x,x)=0$ for all $x\in V$
  \item $b(x,y)>0$ iff $(x,y)\in \vec{E}$.
    \end{itemize}

  In addition, we consider a measure on $V$ given by a positive function 
$$m:V\to \left]  0,\infty\right) .$$
The weighted graph is \textit{symmetric} if  for all $x,y\in V$, $b(x,y)=b(y,x)$,
  as a consequence $(x,y)\in \vec{E}\Rightarrow (y,x)\in\vec{E}$.\\
  The graph is called  \textit{simple} if the weights $m$ and $b$ are constant equal to $1$ on $V$ and $\vec{E}$ respectively.\\

On a non symmetric graph we have three notions of connectedness, see \cite{AT}, \cite{Ma}.
\begin{defini}
\begin{itemize}
\item The set of undirected edges is defined by 
$$E=\left\{\{x,y\},~(x,y)\in\vec{E}\text{ or } (y,x)\in\vec{E}\right\}.$$
\item  A chain from the vertex $x$ to the vertex $y$ in $G$ is a finite set of undirected edges $\{x_1,y_1\};~\{x_2,y_2\};..;\{x_n,y_n\},~n\geq 1$ 
$$x_1=x,~y_n=y \text{ and } x_i=y_{i-1}~~\forall~2\leq i\leq n.$$
\item A path from the vertex $x$ to the vertex $y$ in $G$ is a finite set of directed edges $(x_1,y_1);~(x_2,y_2);..;(x_n,y_n),~n\geq 1$ such that
$$x_1=x,~y_n=y \text{ and } x_i=y_{i-1}~~\forall~2\leq i\leq n.$$
\item $G$ is called weakly connected if two vertices are always related by a chain.
\item $G$ is called connected if two vertices are always related by a path.
\item $G$ is called strongly connected if there is for all vertices $x,y$ a path from $x$ to $y$ and one from $y$ to $x$.
 \end{itemize}
\end{defini}
We assume in the following that the graph under consideration is \textbf{connected}, \textbf{locally finite} and satisfies:
\begin{equation}\label{cco}
\forall x\in V,~\exists y\in V;~(x,y)\in\vec{E}.
\end{equation}

\subsection{Functional spaces}
Let us introduce the following spaces associated to the graph $G$:\\
The space of functions on the graph $G$ is considered as the space of complex
functions on V and is denoted by
$$\mathcal{C}(V)=\{f:V\to \mathbb{C} \}$$
 and $\mathcal{C}_c(V)$ is its subset of finite supported functions.\\
  We consider for a measure $m$, the space
$$\ell^2(V,m)=\{f\in \mathcal{C}(V), ~~\sum_{x\in V}m(x)|f(x)|^2<\infty\}$$
It is a Hilbert space when equipped with the scalar product given by
$$(f,g)_m=\sum_{x\in V}m(x)f(x)\overline{g(x)}.$$
The associated norm is given by:
$$\|f\|_m=\sqrt{(f,f)_m}.$$
\subsection{Laplacians on a directed graph}

For a locally finite, without loops, connected graph and satisfying \eqref{cco}, we introduce the combinatorial Laplacian $\Delta$
defined on $\mathcal{C}_c(V)$ by:
$$\Delta f(x)=\frac{1}{m(x)}\sum_{y\in V}b(x,y)\left( f(x)-f(y)\right) .$$
\textbf{ Dirichlet operator}:\\ Let $U \subset V$, the Dirichlet Laplacian $\Delta^{D}_U$ defined on $\mathcal{C}_c(U)$ is the restriction of $\Delta$ on $U$.\\

In the sequel, we introduce the \textbf{Assumption} $(\beta)$ already used in \cite{Mar}, \cite{Ma} and \cite{M}, which is like Kirchhoff's law in electric networks. \\

  \textbf{Assumption $(\beta)$}: for all $x\in V$,~ $\beta^+(x)=\beta^-(x)$\\
  where $$\beta^+(x)=\sum_{y\in V }b(x,y)\text{ and }\beta^-(x)=\sum_{y\in V}b(y,x).$$
  
In the sequel of this work, we suppose that the Assumption $(\beta)$ is satisfied.\\

 \begin{prop}
If the graph $G$ satisfies the Assumption $(\beta)$, the formal adjoint $\Delta'$ of the operator $\Delta$ is defined on $\mathcal{C}_c(V)$ by:
$$\Delta' f(x)=\frac{1}{m(x)} \sum_{y\in V}b(y,x)\left(f(x)-f(y)\right) .$$
\end{prop}
\begin{demo}
The following calculation gives for all $f,g \in \mathcal{C}_c(V)$:
 \begin{align*}
(\Delta f,g)_m=&\sum_{(x,y)\in \vec{E}}b(x,y)\big( f(x)-f(y)\big) \overline{g(x)}\\
=&\sum_{x\in V}f(x)\overline{g(x)}\sum_{y\in V}b(y,x)-\sum_{(y,x)\in \vec{E}}b(y,x)f(x) \overline{g(y)}\\
=&\sum_{(y,x)\in \vec{E}}b(y,x)f(x)\left(\overline{g(x)}- \overline{g(y)}\right) \\
=&(f,\Delta'g)_m.
\end{align*}
\end{demo}
\begin{rem}
As a consequence, the adjoint operator $\Delta^*$ of $\Delta$ is defined by:
$$D(\Delta^*)=\left\lbrace f\in\ell^2(V,m),~\frac{1}{m(x)} \sum_{y\in V}b(y,x)\left(f(x)-f(y)\right)\in\ell^2(V,m)\right\rbrace. $$
\end{rem}
An explicit Green's formula associated to the non self-adjoint Laplacian $\Delta$ is established in \cite{Ma}:
\begin{lemm}(Green's Formula)\label{grn}
Let $f$ and $g$ be two functions of $\mathcal{C}_c(V)$. Then under the Assumption $(\beta)$ we have
$$(\Delta f,g)_m+(\Delta' g,f)_m=\sum_{(x,y)\in \vec{E}}b(x,y)\big( f(x)-f(y)\big) \big( \overline{g(x)-g(y)}\big).$$
\end{lemm}

\section{ Sectoriality of the Laplacian}\label{2}
We refer to \cite{Kts}, \cite{Ev} for the basic definitions and results concerning
sectorial operators.
\begin{defi} 
 The numerical range of an operator $A$ with domain $D(A)$, denoted by  $W(A)$, is the non-empty set
$$W(A)=\{(A f,f),~~f\in D(A),~\parallel f\parallel=1 \}.$$
\end{defi}

\begin{defi}
Let $\mathcal{H}$ be a Hilbert space, an operator $A: D(A)\to H$ is said to be sectorial if $W(A)$ lies in a sector
$$S_{a,\theta}=\{z\in\mathbb{C},~\mid arg(z-a)\mid \leq \theta \}$$
for some $a\in \mathbb{R}$ and $\theta\in \left[  0,\frac{\pi}{2}\right) $.
\end{defi}
We introduce an assumption on $G$ that allows us to study the sectoriality of $\Delta$. \\ 

\textbf{Assumption $(\gamma)$}:\\
\begin{equation*}
\exists~M>0,~\forall~x\in V,~\sum_{y\in V}\mid b(x,y)-b(y,x)\mid \leq M m(x)
\end{equation*}
where $M$ is a positive constant.\\

\begin{figure}[ht]
\vspace{0.05cm}
\begin{center}
\includegraphics*[height=8cm,width=11cm]{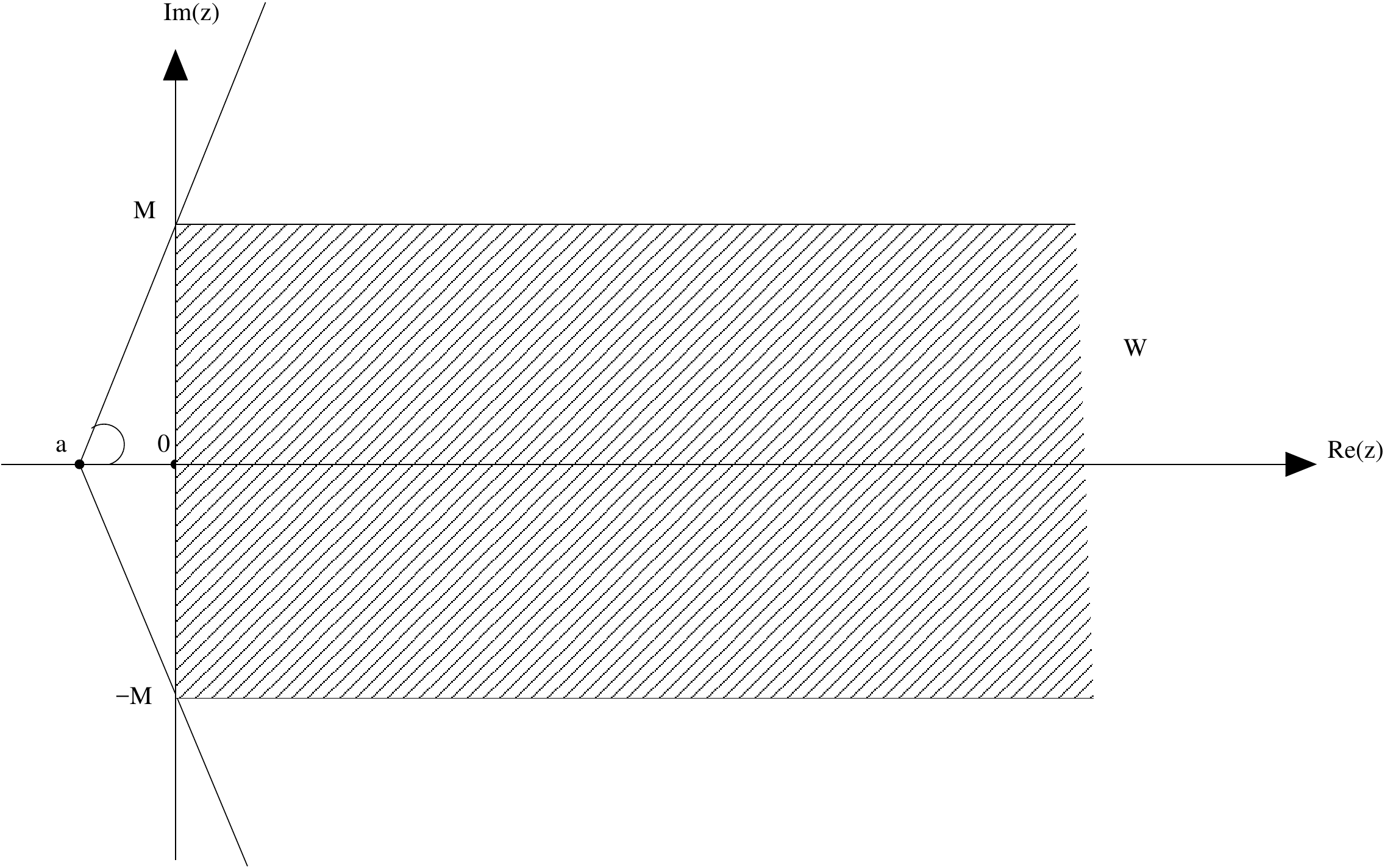}
\end{center}
 \caption{$W$ in the complex plane.}
 \label{fig1}
\end{figure}

Consider the following operator defined on $\mathcal{C}_c(V)$ by :
 $$(\Delta-\Delta')f(x)=\frac{1}{m(x)}\sum_{y\in V}\big(b(x,y)-b(y,x)\big)\big(f(x)-f(y)\big).$$
\begin{prop}\label{boun} If the Assumptions $(\beta)$ and $(\gamma)$ are satisfied, then the operator $(\Delta-\Delta')$, which is defined on $\mathcal{C}_c(V)$, extends to a unique bounded operator on $\ell^2(V,m)$.
\end{prop}
\begin{demo}
Let $f,g\in \mathcal{C}_c(V)$, by the Assumption $(\beta)$, we have
\begin{align*}
(\Delta f,g)-(\Delta'f,g)&=\sum_{x\in V}\beta^+(x)f(x)\overline{g(x)}-
\sum_{x\in V}\beta^-(x)f(x)\overline{g(x)}\\
&+\sum_{x\in V}\overline{g(x)}\sum_{y\in V}\big(b(x,y)-b(y,x)\big)f(y)\\
&=\sum_{x\in V}\overline{g(x)}\sum_{y\in V}\big(b(x,y)-b(y,x)\big)f(y).
\end{align*}
The Cauchy-Schwarz inequality give
\begin{align*}
\big|\big((\Delta-\Delta')f,g\big)\big|&\leq \sum_{x\in V}|\overline{g(x)}|\sum_{y\in V}\big|b(y,x)-b(x,y)\big||f(y)|\\
&\leq \sum_{x\in V}|\overline{g(x)}|\Big(\sum_{y\in V}\big|b(y,x)-b(x,y)\big|\Big)^\frac{1}{2}\\
&~~~~~~~~~~~~~~~~\Big(\sum_{y\in V}\big|b(y,x)-b(x,y)\big||f(y)|^2\Big)^\frac{1}{2}
\end{align*}
now, we apply the Assumption $(\gamma)$
\begin{align*}
\big|\big((\Delta-\Delta')f,g\big)\big|&\leq \Big(M\sum_{x\in V}m(x)|\overline{g(x)}|^2\Big)^\frac{1}{2}\Big(\sum_{x\in V}\sum_{y\in V}\big|b(y,x)-b(x,y)\big||f(y)|^2\Big)^\frac{1}{2}\\
&\leq \Big(M\sum_{x\in V}m(x)|\overline{g(x)}|^2\Big)^\frac{1}{2}\Big(\sum_{y\in V}\sum_{x\in V}\big|b(y,x)-b(x,y)\big||f(y)|^2\Big)^\frac{1}{2}\\
&\leq \Big(M\sum_{y\in V}m(x)|\overline{g(x)}|^2\Big)^\frac{1}{2}\Big(\sum_{y\in V}|f(y)|^2\sum_{x\in V}\big|b(y,x)-b(x,y)\big|\Big)^\frac{1}{2}\\
&\leq \Big(M\sum_{x\in V}m(x)|\overline{g(x)}|^2\Big)^\frac{1}{2}\Big(M\sum_{y\in V}m(y)|f(y)|^2\Big)^\frac{1}{2}.
\end{align*}
Thus $\|(\Delta-\Delta')f\|_m\leq M\|f\|_m$, $\forall f\in\mathcal{C}_c(V)$. Then we conclude by the Theorem of Hahn-Banach.
\end{demo}
\begin{prop}
If the Assumptions $(\beta)$ and $(\gamma)$ are satisfied, then the Laplacian $\Delta$ is sectorial.
\end{prop}
\begin{demo}
From the Green's formula, we have
\begin{align*}
2\mathcal{R}e(\Delta f,f)&=(\Delta f,f)+(\Delta' f,f)\\
&=\sum_{(x,y)\in \vec{E}}b(x,y)\mid f(x)-f(y)\mid^2\\
&\geq 0
\end{align*}

By Proposition \ref{boun}, $\mathcal{I}m(\Delta f,f)$ is bounded by $\dfrac{M}{2}$. Hence the real part of the numerical range is positive and its imaginary part is bounded. Then $\Delta$ is sectorial where $a$ is any point in the open half line of the negative real part.
 See Figure \ref{fig1}.
\end{demo}

Recall that a sectorial operator is closable as an operator if it is densely defined, see Theorem V-3.4 of \cite{Kts}.
\begin{coro}
If the Assumptions $(\beta)$ and $(\gamma)$ are satisfied, then the Laplacian $\Delta$ is closable.
\end{coro}

\section{Lack of the essential spectrum for $\overline{\Delta}$} \label{sym}
In the study of the non self-adjoint problem, it is natural to hope to use the extensive literature that already exists for self-adjoint problems. Hence, we introduce some spectral properties of the symmetric operator $H=\dfrac{1}{2}(\Delta+\Delta')$ as the essential self-adjointness character. We present a comparison theorem concerning the essential spectra of $\overline{\Delta}$ and $S=\overline{H}$. We investigate when the emptiness of the essential spectrum of $S$ implies the emptiness of the essential spectrum of $\overline{\Delta}$.

\subsection{Essential self-adjointness}
The question of the essential self-adjointness of a symmetric Laplacian of graphs is of most importance for the study of infinite graphs and it received many answers. See for instance \cite{kel} and the references inside. We refer here to \cite{To10} who ensure essential self-adjointness for the Laplacian on a symmetric case. Recently C. Ann\'e and N. Torki \cite{AT} proposed a geometric condition who ensure essential self-adjointness for the Gauss-Bonnet operator and for the the Laplacian.
\begin{defi}
A graph $G$ is called with bounded degree if there exists an integer $N$ such that for any $x\in V$
 we have:
 $$\#\{y\in V,~\{x,y\}\in E\}<N.$$
\end{defi}
\begin{defi}
We define the b-weighted distance on $G$, which we denote by $\delta_b$ by :
 $$\delta_b(x,y)= \min_{\gamma\in\Gamma_{x,y}}L(\gamma)$$
where $\Gamma_{x,y}$ is the set of all chains $\gamma$: $\{x_1,x_2\};~\{x_2,x_3\};..;\{x_{n-1},x_n\},~n\geq 1$, linking the vertex $x$ to the vertex $y$; and $L(\gamma)=\displaystyle{\sum_{1\leq i\leq n}\frac{\sqrt{m(x_i)m(x_{i+1})}}{\sqrt{b(x_i,x_{i+1})+b(x_{i+1},x_i)}}}$ the length of $\gamma$.
\end{defi}
Applying the Theorem 6.2 of N. Torki-Hamza \cite{To10} for the symmetric operator $H$ , we find the following result.

\begin{theo}\label{na}
Let $G$ be an infinite connected graph with bounded degree and which is weighted by $m$ on $V$ and $b$ on $\vec{E}$.  We assume that the metric associated to the distance $\delta_b$ is complete. Then the operator $H$ is essentially self-adjoint.
\end{theo}
C. Ann\'e and N. Torki \cite{AT} gave
a general condition on the graph by defining the notion of $\chi$-completeness for the essential self-adjointness of the Gauss-Bonnet operator and the Laplacian.
This condition covers many situations that have been already studied, as \cite{To10} and \cite{kel}.\\
For directed graphs, we propose the following definition of $\chi$-completeness, introduced in \cite{AT} for symmetric graphs.
\begin{defi}
The directed graph $G$ is $\chi$-complete if there exists an increasing sequence of 
finite sets $(B_n)_{n\in\mathbb{N}}$ such that $V =\cup B_n$ and there exist related functions $\chi_n$ satisfying
the following three conditions
\begin{itemize}
\item  $\chi_n\in \mathcal{C}_c(V), ~0 \leq \chi_n \leq 1$
\item $v\in B_n ~\Rightarrow~\chi_n(v) = 1$
\item $\exists C> 0, ~\forall n \in\mathbb{N}, ~\forall x \in V$, $$\dfrac{1}{m(x)}\sum_{(x,y)\in\vec{E}}(b(x,y)+b(y,x))|\chi_n(x)-\chi_n(y)|^2\leq C.$$
\end{itemize}
\end{defi}
The following Theorem is inspired by the Corollary 14 of \cite{AT}.
\begin{theo}\label{CN}
Let $G$ be a connected, locally finite graph. If $G$ is $\chi$-complete, then the
operator $H$ is essentially self-adjoint.
\end{theo}
We introduce an assumption on $G$ which allows to study the self-adjointness of the symmetrized operator $S$.
\begin{defi}
 We define a Laplacian $S$ on $D(\overline{\Delta})\cap D(\Delta^*)$, as the sum of the two non self-adjoint Laplacians $\overline{\Delta}$ and $\Delta^{*}$, given by:
\begin{align*}
Sf(x)&=\dfrac{1}{2}(\overline{\Delta}+\Delta^{*})f(x)\\
&=\frac{1}{m(x)}\sum_{y\in V}\frac{b(x,y)+b(y,x)}{2}\big(f(x)-f(y)\big).
\end{align*}
\end{defi}
\begin{prop}
 The operator $S$, defined above, is a symmetric extension of $H$.
\end{prop}
\begin{demo}
we have $\forall f\in\mathcal{C}_c(V),~~Sf=Hf.$\\
We know that in general,
$$(\overline{\Delta})^*+(\Delta^{*})^*\subset \left(\overline{\Delta})+\Delta^*\right)^*.$$
 $S$ is symmetric because $S^*$ is an extension of $S$. In fact,
\begin{align*}
S&=\dfrac{1}{2}(\overline{\Delta}+\Delta^*)\\
&=\dfrac{1}{2}\left( (\Delta^*)^*+(\overline{\Delta})^*\right) \\
&\subset \dfrac{1}{2}(\Delta^*+\overline{\Delta})^*\\
&=S^*.
\end{align*}
\end{demo}
\begin{prop}
If the graph $G$ satisfies the Assumption $(\beta)$  and $(\gamma)$, then
$$D(\overline{\Delta})\subset D(\Delta^{*})$$
and $\left( S,D(\overline{\Delta})\right) $ is a closed operator.
\end{prop}
\begin{demo}
By the Assumption $(\gamma)$, we have $\Delta-\Delta'$ is extended to a unique bounded operator $B= \overline{\Delta-\Delta'}$ on $\ell^2(V,m)$, (Proposition \ref{boun}). Therefore, as
$\Delta=\left( \Delta-\Delta'\right) +\Delta'$, we have $\overline{\Delta}=B+\overline{\Delta'}$. It follows that 
$$D(\overline{\Delta})\subset D(\overline{\Delta'})\subset D(\Delta^*).$$
For the closeness of $S$, it sufficient to see that $S=\overline{\Delta}-\dfrac{1}{2}B$.
\end{demo}
\begin{rem}
For a linear bounded operator $T$, it is obvious that $\mathcal{R}e(T)=\dfrac{1}{2}\big(T+T^{*}\big)$, but this is not true in general. In our situation the Assumption $(\gamma)$ implies that $S=\mathcal{R}e(\overline{\Delta})$.
\end{rem}

In the sequel of this work, we assume that the essential self-adjointness condition of $H$ is satisfied.
\begin{prop}
Let $G$ be a graph satisfying the Assumptions $(\beta)$  and $(\gamma)$. If $H$ is essentially self-adjoint (for example $G$ satisfied the Hypothesis of Theorems \ref{na} or \ref{CN}. Then $\left( S,D(\overline{\Delta})\right)$ is a self-adjoint operator.
\end{prop}
\begin{demo}
 $S$ is a symmetric closed extension of $H$.
\end{demo}
\subsection{Essential spectrum of the Laplacian}
Recently, many efforts have been made to study when the spectrum of the graph Laplacian is
discrete or the essential spectrum is empty, \cite{Woj}, \cite{Klr}, \cite {Kl}, \cite{Ma}.  \\

Recall the definition of the essential spectrum stated in \cite{Kts}.

\begin{defi}\label{deff}
 The essential spectrum $\sigma_{ess}(A)$ of a closed operator $A$ is the set of all complex numbers $\lambda$ for which the range $R(A-\lambda)$ is not closed or $\dim\ker(A-\lambda)=\infty$.
\end{defi}

 For a non-compact manifold with finite volume $M$, the lower bound of the essential spectrum of a self-adjoint Laplacian $A$ on $M$, $\lambda_1^{ess}(A)$ admits the following characterization given in \cite{R.B}:
$$\lambda_1^{ess}(A)=\lim_{K\to M}~\lambda_1(A^{D}_{K^c}).$$
where $K$ runs over an increasing set of compact subdomains of $M$ such that $\cup K=M$. 

Then in \cite{Klr}, Keller stated a comparable proposition on graphs which allows us to determine the essential spectrum of a self-adjoint  operator via its restriction on the complement of larger and larger sets.
\begin{prop}\label{klr}
Let $G=(V,\vec{E})$ be an infinite graph, $A$ is a self adjoint operator bounded from below, with domain $D(A)$ such that $\mathcal{C}_c(V) \subseteq D(A)\subseteq \ell^2(V,m)$ then:
$$\inf \sigma_{ess}(A)=\lim_{K,~\text{finite}}\inf_{f\in \mathcal{C}_c(V) \atop supp f \subseteq K^c}\frac{(Af,f)}{(f,f)}=\lim_{K}\inf \sigma(A_K)$$
$$\sup \sigma_{ess}(A)\leq \lim_{K,~\text{finite}}\sup_{f\in \mathcal{C}_c(V) \atop supp f \subseteq K^c}\frac{(Af,f)}{(f,f)}=\lim_{K}\sup \sigma(A_K)$$
where $A_K$ is the restriction of $A$ to $K^c$ with Dirichlet condition.\\
\end{prop}

Let us define the following numbers:

\begin{equation*}
\nu(A)=\inf\{\mathcal{R}e \lambda:~~\lambda\in W(A)\}.
\end{equation*}
\begin{equation*}
\eta^{ess}(A)=\inf\{\mathcal{R}e \lambda:~~\lambda \in \sigma_{ess}(T)\}.
\end{equation*}

If $A$ is a bounded operator, then the spectrum is always a subset of the closure of the numerical range but this is not true in general. Certainly, the essential spectrum of a closed operator is a subset of the closure of the numerical range. Hence, we have for the Laplacian $\overline{\Delta}$:
\begin{equation}\label{rang}
\eta^{ess}(\overline{\Delta})\geq \nu(\overline{\Delta})
\end{equation}
The following Theorem follows from part (IV), Theorem 1.11 of \cite{Ev}. From the sectoriality of $\Delta$, we will be able to  compare the essential spectrum of $\Delta$ and the essential spectrum of its real part.
\begin{theo}
If the Assumptions $(\beta)$ and $(\gamma)$ are satisfied and $H$ is essentially self-adjoint, then
$$\eta^{ess}(\overline{\Delta})\geq \inf \sigma_{ess}(S).$$
\end{theo}

 We recall the definitions of the Cheeger constants on $\Omega\subset V$:
 \begin{equation*}
h(\Omega)~ = \inf_{U\subset \Omega\atop finite}\frac{b(\partial_E U )}{m(U)}
\end{equation*}
and
\begin{equation*}
\tilde{h}(\Omega)~ = \inf_{U\subset \Omega\atop finite}\frac{b(\partial_E U )}{\beta^+(U)}
\end{equation*}
where $\partial_E \Omega=\Big\{(x,y)\in \vec{E}:~ (x\in \Omega,~y\in \Omega^{c} )~~or ~~(x\in \Omega^{c},~y\in \Omega)\Big\}.$\\

  We provide the Cheeger inequality at infinity on a filtration of graph $G$.
 \begin{defi}
A graph $H=(V_H,\vec{E}_H)$ is called a subgraph of $G=(V_G,\vec{E}_G)$ if $V_H\subset V_G$
and $\vec{E}_H=\big\{(x,y)\in V_H\times V_H~~\big\}\cap \vec{E}_G $.
\end{defi}
\begin{defi}
 A filtration of $G=(V,\vec{E})$ is a sequence of finite connected subgraphs $\{G_n=(V_n,\vec{E}_n),~~n\in \mathbb{N} \}$ such that $V_n\subset V_{n+1}$ and:
 $$\displaystyle\cup_{n\geq 1} V_n=V.$$
\end{defi}
Let $(G,b)$ be a connected, weighted infinite graph, $\{G_n,~n\in \mathbb{N}\}$ a filtration of $G$.
Let us denote
$$ M_{V_n^c}=\sup \left\{\frac{\beta^+(x)}{m(x)},~~x\in V_n^c\right\}$$

 The isoperimetric constant at infinity is given by
$$h_\infty=\lim_{n\rightarrow \infty}h(V^{c}_n).$$
 We refer to Theorem 3.7 of \cite{Ma} for the Cheeger estimation of the numerical range of $\Delta^D_\Omega$.
 \begin{theo}\label{chee}
Let $\Omega\subset V$, the bottom of the real part of $W(\Delta^{D}_\Omega)$ satisfies the following inequality:
\begin{equation}\label{co}
\frac{{h}^{2}(\Omega)}{8}~\leq~M_\Omega\nu(\Delta^{D}_\Omega)~\leq ~\frac{1}{2}M_\Omega h(\Omega).
\end{equation}
\end{theo}

In the sequel of this work we suppose that the \textbf{Assumptions} $(\beta)$ and $(\gamma)$ are \textbf{satisfied} and $H$ is \textbf{essentially self-adjoint}.
 
\begin{lemm}\label{S} Let $\Omega$ a subset of $V$. Then the bottom of the spectrum $\lambda_1(S^{D}_\Omega)$ of $S^{D}_\Omega$ satisfies
$$\lambda_1(S^{D}_\Omega)=\nu(\Delta^{D}_\Omega).$$
\end{lemm}
\begin{demo}
The bottom of the spectrum $\lambda_1$ of the Dirichlet Laplacian on $\Omega\subset V$, $S^D_\Omega$ admits the variational definition:
 \begin{align*}
 \lambda_1(S^{D}_\Omega)=& \lambda_1(\overline{H}^{D}_\Omega)\\
:=&\inf_{f\in \mathcal{C}_c(\Omega)\atop \|f\|_m=1}\big(H^D_\Omega f, f\big)_m \\
=&\inf_{f\in\mathcal{C}_c(\Omega)\atop \|f\|_m=1}\mathcal{R}e(\Delta^D_\Omega f,f)\\
=&\nu(\Delta^{D}_\Omega).
 \end{align*}
  \end{demo}
  
We now have all the tools to provide Cheeger's Theorem associated with the self-adjoint operator $S$.
\begin{prop}\label{ba}
Let $\Omega$ be a subset of $V$. Then
\begin{equation}
M_{\Omega}\lambda_{1}(S^{D}_\Omega)\geq\frac{h^2(\Omega)}{8}.
\end{equation}
\end{prop}
\begin{demo}
It is a simple deduction of the Lemma \ref{S}, $$\lambda_{1}(S^{D}_{\Omega})=\nu(\Delta^{D}_{\Omega}).$$
We use the inequality (\ref{co}) to conclude.
\end{demo}

In order to study the essential spectrum of $S$ and the dependence on the geometry at infinity, we start with the following fundamental Lemma.
\begin{lemm}\label{lim}
Let $\{G_,~n\in\mathbb{N}\}$ be a filtration of $G$, we have:
$$\lim_{n\rightarrow\infty }\big(\lim_{k\rightarrow\infty\atop k\geq n+1}\lambda_1(S^{D}_{ G_k\setminus G_n})\big)=\lambda^{ess}_1(S).$$
\end{lemm}
\begin{demo}
Denote by $l=\displaystyle{\lim_{n\rightarrow\infty}\big(\lim_{k\rightarrow\infty\atop k\geq n+1}\lambda_1(S^{D}_{ G_k\setminus G_n})\big)}$, these limits exist because $(G_k\setminus G_n)_{k\geq n+1}$ and $(G_n^c)_n$ are monotones. For each $n\in\mathbb{N}$,  $(G_k\setminus G_n)_{k\geq n+1}$ is a sequence of finite subgraphs whose union is equal to $G^{c}_n$ hence from \cite{akch}, Theorem 2.3.6 we have :
$$\lambda_1(S^D_{G^{c}_n})= \lim_{k\rightarrow\infty}\big(\lambda_1(S^{D}_{ G_k\setminus G_n})\big)$$
and then by the application of Proposition \ref{klr}, we obtain
 $$l=\lambda^{ess}_1(S).$$
\end{demo}

H. Donnelly and P. Li  \cite{Don} showed that the essential spectrum of the Laplacian depends on the geometry at infinity and it is empty on a rapidly curving manifold. M. Keller \cite{Klr} gave a sufficient condition  for the discreteness of the spectrum. This condition is the positivity of a Cheeger constant at infinity $\tilde{h}_\infty$ on a simple  rapidly branching graph. For the weighted graph, we propose a condition which implies the emptiness of the essential spectrum and this condition can be satisfied even if $\tilde{h}_\infty=0$. It generalizes  the work of \cite{Ma} for heavy end graphs where was assumed that $\tilde{h}_\infty>0$.
\begin{theo}\label{Ma}
Let $\{G_n,~n\in\mathbb{N}\}$ be a filtration of $G$, if there exists a sequence $(c_n)_n$ such that for all $k\geq n+1$
\begin{equation}\label{Abs}
\frac{h^{2}(G_k \setminus G_n)}{8 M_{G_k \setminus G_n}}\geq c_n
\text{ and } \lim_{n\rightarrow\infty}c_n=\infty
\end{equation}
 then $\sigma_{ess}(S)$ is empty.
\end{theo}
\begin{demo}
By the Proposition \ref{ba}, we have for all $k\geq n+1$:
$$\lambda_{1}(S^{D}_{G_k\setminus G_n})\geq c_n$$
therefore
$$\lim_{n\rightarrow\infty}\big(\lim_{k\rightarrow\infty}\lambda_1(S^{D}_{ G_k\setminus G_n})\big)\geq \lim_{n\rightarrow\infty} c_n$$

thus from the Lemma \ref{lim} we obtain the result.
\end{demo}

The aim of the sequel of this section is to describe the relationship between the essential spectrum of $\Delta$ and the essential spectrum of its real part, see \cite{Lws}.\\

 Recall now the definition considered by Lewis:
The essential spectrum of a closed operator densely defined $T$  is the set of all complex number $ \lambda $ for which $ T-\lambda I $ has a singular sequence. This definition is equivalent to that given in the Definition \ref{deff}  (see Theorem 1.6 \cite{Ev}).
\begin{theo}[\textbf{Comparison Theorem of Lewis}]\cite{Lws}
Let $T$ be a closed linear operator in the Hilbert space $\mathcal{H}$ with dense domain $D(T)$. Let $A$ be a self-adjoint operator in $\mathcal{H}$ bounded from below and with $D(T)\subset D(A)$. If
$$\mathcal{R}e(Tu,u)\geq (Au,u),~~~~\forall u\in D(T)$$
then
$$\sigma_{ess}(T)\subseteq \{\lambda \in \mathbb{C}:~~\mathcal{R}e(\lambda) \geq \inf\sigma_{ess}(A)\}.$$
If $\sigma_{ess}(A)=\emptyset$ then $\sigma_{ess}(T)=\emptyset$.
\end{theo}
As a consequence we have :
\begin{prop}\label{B}
$$\sigma_{ess}(\overline{\Delta})\subseteq \{\lambda \in \mathbb{C}:~~\mathcal{R}e(\lambda) \geq \inf\sigma_{ess}(S)\}$$
and if $\sigma_{ess}(S)=\emptyset$, then $\sigma_{ess}(\overline{\Delta})=\emptyset$.
\end{prop}
\begin{theo}
Let $G$ be a graph which satisfies the Assumptions $(\beta)$ and $(\gamma)$, that $H$ is essentially self-adjoint and also the Hypothesis (\ref{Abs}). Then $\sigma_{ess}(\overline{\Delta})=\emptyset$.
\end{theo}
\begin{demo}
From the Theorem \ref{Ma} and the Proposition \ref{B}, we deduce the outcome.
\end{demo}
\section{Example}\label{EX}
We propose an example of graph with $\tilde{h}_\infty=0$ and a totally discrete spectrum of $S$, then with the lack of essential spectrum for $\overline{\Delta}$ also.
Let us consider the graph $\mathbb{Z}$ with $G_n=\{-n,...,-1,0,1,2,...,n\}$ as a filtration. The weights on the graph are  $m(l)=1$, $b(l,l+1)=\dfrac{(|l|^{3}+1)}{2}+\dfrac{1}{4}$ and $b(l+1,l)=\dfrac{(|l|^{3}+1)}{2}-\dfrac{1}{4}$, $\forall ~l\in \mathbb{Z}$.
\begin{figure}[ht]
\vspace{0.07cm}
\begin{center}
\includegraphics*[height=0.6cm,width=12cm]{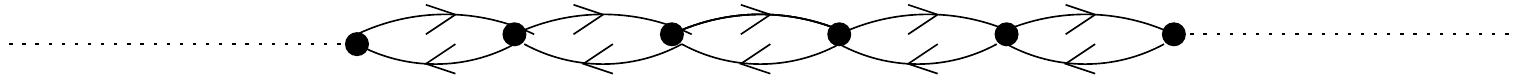}
\end{center}
 \caption{The graph $\mathbb{Z}$}
\end{figure}

Independently of the metric on $G$, we cite the following Theorem inspired by \cite{To10}.
\begin{theo}\label{esnt}
If the weight $m$ is  constant  on $V$ then the operator $H$ is essentially self-adjoint.
\end{theo}
\begin{coro}
$\overline{\Delta}+\Delta^*$ is a self-adjoint operator in $D(\overline{\Delta})$.
\end{coro}
\begin{demo}
 It is clear that the Assumptions $(\beta)$ and $(\gamma)$ are satisfied on $\mathbb{Z}$. From Theorem \ref{esnt} $H$ is essentially self-adjoint on $D(\overline{\Delta})$.  \\
\end{demo}

 Now, we fix $n$ and let $k>n$, we have for $k=2n$:
\begin{align*}
  \tilde{h}(G^c_n)&\leq\dfrac{b(\partial_E(G_k \setminus G_n))}{\beta^+(G_k \setminus G_n)}\\
  &\leq   \dfrac{n^{3}+(k+1)^{3}}{\sum_{t=n}^{k}\frac{1}{2}\left(t^{3}+(t+1)^{3}\right)}\\
&\leq  \dfrac{2\big(n^{3}+(2n+1)^{3}\big)}{(n+1)\big(n^{3}+(n+1)^3\big)}.
\end{align*}
Therefore
 $$ \tilde{h}_\infty =\lim_{n\rightarrow \infty}\tilde{h}(G^c_n)\leq \lim_{n\rightarrow \infty}\dfrac{2\big(n^{3}+(2n+1)^{3}\big)}{(n+1)\big(n^{3}+(n+1)^3\big)}=\lim_{n\rightarrow \infty} \dfrac{18n^{3}}{2n^{4}}=0.$$\\
 
In addition, $\sigma_{ess}(S)=\emptyset$ on $\mathbb{Z}$, in fact:\\
Let $f$ be a function with finite support, we have
\begin{align*}
\|f\|^2_m=&\sum_{k\in \mathbb{Z}}|f(k)|^2\\
=&\sum_{k\in \mathbb{Z}}(k+1-k)|f(k)|^2\\
=&\sum_{k\in \mathbb{Z}}(k+1)|f(k)|^2-\sum_{k\in \mathbb{Z}}k|f(k)|^2\\
=&\sum_{k\in \mathbb{Z}}k(|f(k-1)|^2-|f(k)|^2)\\
=&\sum_{k\in \mathbb{Z}}k(|f(k-1)|-|f(k)|)(|f(k-1)|+|f(k)|).
\end{align*}
By Cauchy-Schwarz inequality, we obtain
\begin{align*}
\|f\|^2_m\leq &\Big(\sum_{k\in \mathbb{Z}}k^2|f(k-1)-f(k)|^2\Big)^\frac{1}{2}\Big(\sum_{k\in \mathbb{Z}}(|f(k-1)|+|f(k)|)^2\Big)^\frac{1}{2}\\
\leq &\Big(\sum_{k\in \mathbb{Z}}k^2|f(k-1)-f(k)|^2\Big)^\frac{1}{2}\Big(4\sum_{k\in \mathbb{Z}}|f(k)|^2\Big)^\frac{1}{2}
\end{align*}
therefore
$$\dfrac{1}{4}\|f\|^2_m \leq \sum_{k\in \mathbb{Z}} k^2|f(k-1)-f(k)|^2.$$
Now, if $f$ is with support in $V_n^c$, $|k|\geq n$,
$$\dfrac{1}{4}\|f\|^2_m\leq \sum_{|k|\geq n}\dfrac{(|k|^3+1)}{n}|f(k-1)-f(k)|^2.$$
Hence for all non zero function $f$ and with finite support in $V_n^c$, we have
$$\dfrac{n}{8}\leq\dfrac{(S^D_{V_n^c}f,f)_m}{(f,f)_m}=\dfrac{\displaystyle{\sum_{|k|\geq n}\dfrac{|k|^3+1}{2}|f(k-1)-f(k)|^2}}{(f,f)}.$$
This implies that if $n\to \infty$, $\lambda_1(S^D_{V_n^c})\to \infty$. Thus from the Proposition \ref{klr}, we deduce that $\sigma_{ess}(S)=\emptyset$ and therefore the absence of the essential spectrum of $\overline{\Delta}$.\\

\thanks{
\textbf{Acknowledgments}: The author Marwa Balti have benefited to a financial support from the D\'efiMaths program of the Federation of Mathematical Research of the "Pays de Loire" and from the "PHC Utique" program of the French Ministry of Foreign Affairs and Ministry of higher education and research and the Tunisian Ministry of higher education and scientific research in the CMCU project number 13G1501 "Graphes, G\'eom\'etrie et Th\'eorie Spectrale" during her visits to the Laboratory of Mathematics Jean Leray of Nantes (LMJL).
 Also, the three authors would like to thank the Laboratory of Mathematics Jean Leray of Nantes (LMJL) and the research unity (UR/13ES47) of Faculty of Sciences of Bizerta (University of Carthage) for their continuous financial support.}


\begin{thebibliography}{777}
\lhead{{Bibliographie}}\rhead{\textit{Marwa Balti}}
\bibitem[A10]{akch} S. Akkouche.
\newblock {\it Sur la th\'eorie spectrale des op\'erateurs de Schr\"odinnger discrets.}
\newblock Th\`ese de Doctorat, Universit\'e de Bordeaux (2010).
\\

\bibitem[AT15]{AT} C. Ann\'e, N. Torki-Hamza
\newblock {\it The Gauss-Bonnet operator of an infinite graph.}
\newblock Anal. Math. Phys.  \textbf{5}, (2015), 137-159.
\\
\bibitem[AE12]{Ar} W.  Arendt  and  A.F.M.  ter  Elst.
\newblock {\it From  forms  to  semigroups.}
\newblock Spectral theory, mathematical system  theory, evolution equations, differential and difference equations, Oper. Theory Adv. Appl., \textbf{221}, Birkh\"auser/Springer Basel AG, Basel, (2012), 47-69.
\\
\bibitem[Bal17]{Mar} M. Balti.
\newblock {\it Laplaciens non auto-adjoints sur un graphe orient\'e.}
\newblock Th\`ese de Doctorat, Universit\'e de Carthage et Universit\'e de Nantes (2017).
\\
\bibitem[Ba17]{M} M. Balti.
\newblock {\it On the eigenvalues of weighted directed graphs.} 
\newblock Complex analysis and operator theory, \textbf{11}, (2017), 1387-1406.
\\

\bibitem[B17]{Ma} M. Balti.
\newblock {\it Non self-adjoint Laplacians on a directed graph.}  Filomat, \textbf{18}, (2017), 5671-5683.
\\
\bibitem[B84]{R.B} R. Brooks.
\newblock {\it On the spectrum of non-compact manifolds with finite volume.}
\newblock Mathematische Zeitschrift, \textbf{187}, (1984), 425-432.
\\
\bibitem[DL79]{Don} H. Donnelly, P. Li.
\newblock {\it Pure point spectrum and negative curvature for noncompact manifolds.}
\newblock Duke Math J, \textbf{46}, (1979), 497-503.
\\
\bibitem[ELZ83]{Ev} W. D. Evans, R. T. Lewis, A. Zettl.
\newblock {\it Non self-adjoint operators and their essential spectra.
In form local times to global geometry, control and physics, D. Ellworthy, ed.,}
\newblock Differential Equation and Operators, Lecture Notes in Mathematics, \textbf{1032}, (1983), 123-160.
\\
\bibitem[HJMW13]{kel} X. Huang, M. Keller, J. Masamune, R.K. Wojciechowskil.
\newblock {\it A note on self-adjoint
extensions of the Laplacian on weighted graphs.}
\newblock J. Funct. Anal. \textbf{265}, (2013), 1556-1578.
\\
\bibitem[Kat76]{Kts} T. Kato.
\newblock {\it Perturbation theory for linear operators.}
\newblock Springer-Verlag, Berlin, Heidelburg and New York, (1976).
\\
\bibitem[K10]{Klr}  M. Keller.
\newblock {\it The essential spectrum of the Laplacian on rapidly branching tessellations.}
\newblock Mathematische Annalen \textbf{346}, (2010), 51-66.
\\
\bibitem[KMP16]{Kl}  M. Keller, F. M\"unch, AND F. Pogorzelski.
\newblock {\it Geometry and spectrum of rapidly branching graphs.}
\newblock Mathematische Nachrichten \textbf{289}, (2016), 1636-1647 .
\\

\bibitem[Kh13]{Kh} M. Khanalizadeh.
\newblock {\it Sectorial forms and m-sectorial operators.}
\newblock Seminararbetit zum fach funktionalanalysis, Technische Universit\"at Berlin, (2013).
\\
\bibitem[F96]{Fuji} K. Fujiwara.
\newblock {\it The Laplacian on rapidly branching trees.}
\newblock Duke Mathematical Journal.
\textbf{83}, (1996), 191-202.
\\
\bibitem[L79]{Lws} R. T. Lewis.
\newblock {\it Applications of a comparison for quasi-accretive operators in a Hilbert space.}
\newblock In: Everitt W., Sleeman B. (eds) Ordinary and Partial Differential Equations. Lecture Notes in Mathematics. Springer, Berlin,
Heidelberg. \textbf{964}, (1982), 422-434.
\\

\bibitem[T-H10]{To10} N. Torki-Hamza.
 \textit{ Laplaciens de graphes infinis. I: Graphes m\'etriquement complets.}
 Confluentes Math. \textbf{2} (2010), 333-350.\\
Translated to: 
\newblock Essential self-adjointness for combinatorial Schr\"odinger
operators I- Metrically complete graphs.
  arXiv:1201.4644v1.
 \\
 \bibitem[Y74]{Yos} K. Yosida.
\newblock {\it Functional Analysis.}
\newblock Springer-Verlag, Berlin Heidelberg New York, (1974).
\\
 \bibitem[W09]{Woj}   R. K. Wojciechowski.
\newblock {\it Heat kernel and essential spectrum of infinite graphs.}
\newblock Indiana Univ. Math. J. \textbf{58}, (2009), 1419-1441.
\end{thebibliography}
\end{document}